%% file: bcd_pep_long.tex
\documentclass{article}

\usepackage{algorithm, algpseudocode}
\usepackage{amsfonts, amsmath, amsopn, amsthm, epsfig}
  \numberwithin{equation}{section}
\usepackage{graphicx, epstopdf}
\usepackage{hyperref}
\usepackage{subfigure}
\usepackage{appendix}

\usepackage{calc}
\usepackage[top=2.5cm,bottom=2.3cm,inner=2.3cm,outer=2.3cm,bindingoffset=0.5cm,footskip=0.55cm]{geometry}
\newtheorem{theorem}{Theorem}[section]

\newtheorem{lemma}[theorem]{Lemma}

\theoremstyle{remark}
\newtheorem{remark}{Remark}

\newcommand{\ve}[2]{\langle #1 ,  #2 \rangle}

\newenvironment{lemma*}[2][Lemma]{\par\bgroup{\bfseries #1\ #2. }\it\ignorespaces}{\egroup}

\title{A better convergence analysis of the block coordinate descent method for large scale machine learning}

\author{Ziqiang Shi\footnotemark[1] \footnotemark[2], Rujie Liu\footnotemark[1]
}

\input{macros.tex}

\begin{document}

\maketitle

\renewcommand{\thefootnote}{\fnsymbol{footnote}}

\footnotetext[1]{Fujitsu Research \& Development Center, Beijing, China.}
\footnotetext[2]{shiziqiang@cn.fujitsu.com}

\renewcommand{\thefootnote}{\arabic{footnote}}

\begin{abstract}
This paper considers the problems of unconstrained minimization of large scale smooth convex functions having
block-coordinate-wise Lipschitz continuous gradients. The block coordinate descent (BCD) method are among the first optimization schemes suggested for solving such problems~\cite{nesterov2012efficiency}. We obtain a new lower (to our best knowledge the lowest currently)  bound that is $16p^3$ times smaller than the best known on the information-based complexity of BCD method based on an effective technique called Performance Estimation Problem (PEP) proposed by Drori and Teboulle~\cite{drori2012performance} recently for analyzing the performance of first-order black box optimization methods. Numerical test  confirms our analysis.
\end{abstract}
%=============================================

\section{Introduction and problem statement}
\label{sec:introduction}

In this work, we consider the block coordinate descent (BCD) algorithms for solving the large scale problems
of the following form:
\begin{align}
  \min_{x \in \mathbb{R}^D} \,f(x),
  \label{eq:problem}
\end{align}
where $f(x)$ is a smooth convex function (no need to be strongly convex), and it is assumed throughout this work that
\begin{itemize}
\item The gradients of $f(x)$ are
block-coordinate-wise Lipschitz continuous with const $L_i$ $(i=1,\cdots, p)$
\begin{equation}\label{eq:block_coordinate_lipschitz}
\|\nabla_i f(x+\mathbf{U}_i h_i)-\nabla_i f(x+\mathbf{U}_i h_i)\|\leq L_i \|h_i\|,
\end{equation}
where $\mathbf{I}_D=[\mathbf{U}_1,\dots,\mathbf{U}_p]$ is a decomposition of the $D\times D$ identity matrix $\mathbf{I}_D$ into column submatrices
$\mathbf{U}_i\in \mathbb{R}^{D\times D_i}$, and the space $\mathbb{R}^D$ is decomposed into $p$
subspaces: $\mathbb{R}^D=\mathbb{R}^{D_1}\times\dots\times \mathbb{R}^{D_p}$, while $\nabla_i f(x)\in \mathbb{R}^ {D_i}$ is the block of partial derivatives $\nabla_i f(x)=\mathbf{U}_i^\top \nabla f(x)$. We denote the set of functions satisfy this condition as $\mathcal{F}_{\mathbf{L},\mathbf{U}}(\mathbb{R}^D)$, here $\mathbf{L}$ stands for $(L_1,L_2,...,L_p)$, and $\mathbf{U}$ stands for $(\mathbf{U}_1,\dots,\mathbf{U}_p)$.
\item The optimal set $X_*(f):=\arg\min f(x)$ is nonempty, i.e., the Problem~\eqref{eq:problem} is solvable.
\end{itemize}

Block coordinate descent (BCD) methods have recently gained in popularity for solving the Problem~(\ref{eq:problem}) both in theoretical
optimization and in many applications, such as machine learning, signal processing,
communications, and so on. These problems are of very large scale, and
the computational is simple and the cost is very cheap per iteration of BCD methods, yielding computational efficiency.
If moderate accuracy solutions are sufficient for the target applications, BCD methods are often the best option to solve the Problem~(\ref{eq:problem}) in a reasonable time. For convex optimization problems, there exists an extensive
literature on the development and analysis of BCD methods, but most of them focus on the randomized BCD methods~\cite{richtarik2011iteration,nesterov2012efficiency,wright2015coordinate,lu2015on,nutini2015coordinate},
where blocks are randomly chosen in each iteration. In contrast, existing literature on cyclic BCD
methods is rather limited~\cite{beck2013on,Li2016An}, and the later~\cite{Li2016An} is focused on \textbf{strongly} convex functions. In this paper, we focuse on the theoretical
performance analysis of cyclic BCD methods for unconstrained minimization with an objective
function which is known to satisfy the assumptions in the Problem~(\ref{eq:problem})
over the Euclidean space $\mathbb{R}^D$, although the function itself is not known.

We consider finding a minimizer over $\mathbb{R}^D$ of a cost function $f$ belonging to
the set $\mathcal{F}_{\mathbf{L},\mathbf{U}}(\mathbb{R}^D)$.
The class of standard and popular cyclic algorithms of interest generates a sequence of
points $\{x_k^i\in\mathbb{R}^D\;:\;k=0,\cdots,N, \quad i=0,\cdots,p\}$ using the following scheme:

\begin{algorithm}[H]
\caption{The cyclic BCD method}
\label{alg:bcd}

\textbf{Input}: start point $x_0 \in \mathbb{R}^D$.

1: \textbf{repeat for $k=0,1,...,N$}

2: \quad Set $x_k^0=x_k$, and generate recursively
\begin{equation}\label{eq:update_x_ki}
x_k^i=x_k^{i-1}-\frac{1}{L_i}\mathbf{U}_i \nabla_i f(x_k^{i-1})
\end{equation}
\quad \quad for $i=1,...,p$.

3: \quad Update: $x_{k+1}=x_k^p$.

\textbf{Output}: $x_{k+1}$.
\end{algorithm}
The update step at the $k$th iterate $x_k$ performs a gradient step with constant stepsize with respect to a different block of
variables taken in a cyclic order. Evaluating the convergence bound of such BCD algorithms is essential. The sequence $\{x_k\}$ is known to satisfy the bound~\cite{beck2013on}:
\begin{align}
f(x_k)-f(x_*)\leq 4L_{\text{max}}(1+pL^2/L_{\text{min}}^2)R^2(x_0)\frac{1}{k+8/p}
\label{eq:bcd_conv_beck}
\end{align}
for $k\geq 0$, which to our best knowledge is the previously best know analytical bound of cyclic BCD method for
unconstrained smooth convex  minimization. Here $L_{\text{max}}$ and $L_{\text{min}}$ are the maximal and minimal block Lipschitz constants
\begin{equation}\label{eq:lmax}
L_{\text{max}}=\max_{i=1,...,p}\quad L_i \quad \text{and} \quad L_{\text{min}}=\max_{i=1,...,p}\quad L_i,
\end{equation}
$L$ is Lipschitz constant of $\nabla f(x)$, that is
\begin{equation}\label{eq:delta_f_lipschitz}
\|\nabla f(x) - \nabla f(y)\|\leq L\|x - y\|.
\end{equation}
for every $x, y \in \mathbb{R}^D$, and $R(x_0)$ is define by
\begin{equation}\label{eq:rx0}
R(x_0):=\max_{x\in \mathbb{R}^D}\max_{x_*\in X_*(f)}\{\|x-x_*\|:f(x)\leq f(x_0)\}
\end{equation}
same as~\cite{nesterov2012efficiency,beck2013on}.

\begin{figure}[!h]
\subfigure[$p$=2]{
\includegraphics[width=0.48\textwidth]{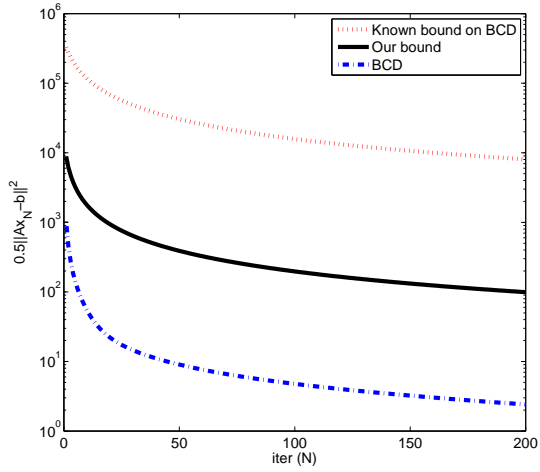}
}
%\hspace{0.1in}
\subfigure[$p$=5]{
\includegraphics[width=0.48\textwidth]{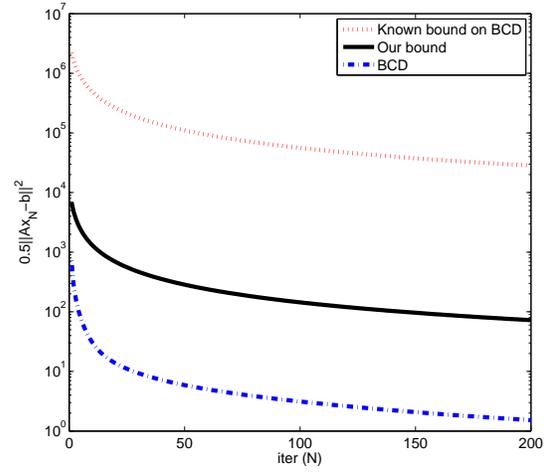}
}
\subfigure[$p$=20]{
\includegraphics[width=0.48\textwidth]{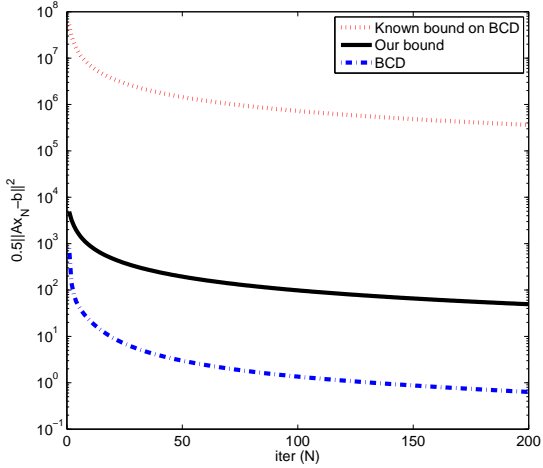}
}
\subfigure[$p$=100]{
\includegraphics[width=0.48\textwidth]{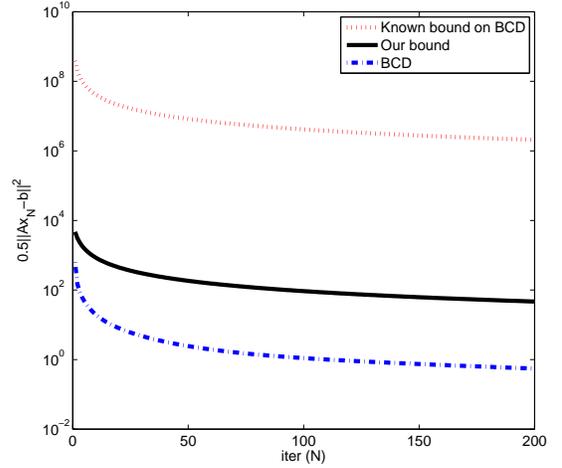}
}
\caption{ The new bounds on the BCD. }
\label{fig:new_bounds}       % Give a unique label
\end{figure}

But in practice, BCD converges much fast. It can be seen in Figure~\ref{fig:new_bounds} in Section~\ref{sec:num_test} tha there is big gap between the currently best known bound and the practice convergence. This work is try to fill this gap.
Recently, Drori and Teboulle~\cite{drori2012performance} considered the Performance Estimation Problem (PEP) approach to
bounding the decrease of a cost function $f$. Following this excellent work, we can formulate the worst case performance bound of the
BCD method over all smooth convex functions $f\in \mathcal{F}_{\mathbf{L},\mathbf{U}}(\mathbb{R}^D)$ as the solution of the following constrained optimization problem:
\begin{align*}
    \begin{aligned}
    \max_{f\in \mathcal{F}_{\mathbf{L},\mathbf{U}}(\mathbb{R}^D)} \max_{ x_0^0, \ldots,x_k^i,\ldots, x_N^p, x_\ast \in\mathbb{R}^D} \quad & f(x_N^p)-f(x_*) \\
    \text{s.t. }\quad
                        & x_k^i=x_k^{i-1}-\frac{1}{L_i}\mathbf{U}_i \nabla_i f(x_k^{i-1}),\quad k=0,\ldots, N,\quad i=1,\ldots, p,\\
                        & x_k^p=x_{k+1}^0,\quad k=0,\ldots, N-1, \\
                        & x_* \in X_*(f).
    \end{aligned}
    \tag{P}\label{P:baseproblem}
\end{align*}

%=============================================
\subsection{Lemmas}
\label{sec:notations}

%Most notations have been introduced in the above, here we mention some facts, which will be used in this work here:
%\begin{equation}\label{eq:summation_of_ux}
%x=\sum_{i=1}^{p} \mathbf{U}_i x^i  \quad \text{and} \quad \sum_{i=1}^{p} \mathbf{U}_i \mathbf{U}_i^\top=\mathbf{I}_D.
%\end{equation}

In the sequel, we often need to estimate from above the differences between two block partial gradients. For that it is convenient to use the following simple lemma:
\begin{lemma}\label{lem:bound_diff_nabla}
Let $f(x)\in \mathcal{F}_{\mathbf{L},\mathbf{U}}(\mathbb{R}^D)$, then we have
\begin{equation}\label{eq:bound_diff_nabla_i}
\frac{1}{2L_i}\|\nabla_i f(y) - \nabla_i f(x)\|^2\leq f(y)-f(x)-\ve{\nabla f(x)}{y-x}
\end{equation}
for every $x, y \in \mathbb{R}^D$.
\end{lemma}

We also need the following lemma (similar but different from Lemma 3.1 of~\cite{drori2012performance}) to simplify a quadratic function of matrix variable into a function of vector variable.
\begin{lemma}\label{lem:quadratic_min}
Let $f(X)=\trace (AX^\top B X + 2 ba^\top X)$ be a quadratic function, where $X\in
\mathbb{R}^{n\times m}$, $A, B\in \mathbb{S}^m$, $a\in \mathbb{R}^n$ and $0\neq b\in
\mathbb{R}^m$. Then
\[
    \inf_{X\in \mathbb{R}^{n\times m}} f(X) = \inf_{\xi\in \mathbb{R}^n} f(\xi b^\top).
\]
\end{lemma}

The proofs of these two lemmas are contained in the appendix for the completeness of this work.

\section{Relaxations of the PEP}

Since Problem~(\ref{P:baseproblem}) involves an unknown function $f$ as a variable, PEP is infinite-dimensional. Nevertheless, it can be
relaxed by using the property of the functions belong to $\mathcal{F}_{\mathbf{L},\mathbf{U}}(\mathbb{R}^D)$.

Let $\{x_0^0,...,x_0^p,...,x_N^p\}$ be the sequence generated by the Algorithm~(\ref{alg:bcd}). Applying (\ref{eq:bound_diff_nabla_i}) to $\nabla_t f(x_m^i)$, $\nabla_t f(x_n^j)$, and $\nabla_t f(x_*)$, and note (\ref{eq:update_x_ki}), we get
\begin{eqnarray*}
\frac{1}{2L_t}\|\nabla_t f(x_m^i) - \nabla_t f(x_n^j)\|^2 &\leq& f(x_m^i)-f(x_n^j)-\ve{\nabla f(x_n^j)}{x_m^i-x_n^j}, \\
\frac{1}{2L_t}\|\nabla_t f(x_m^i) - \nabla_t f(x_*)\|^2 &\leq& f(x_m^i)-f(x_*)-\ve{\nabla f(x_*)}{x_m^i-x_n^j}, \\
\frac{1}{2L_t}\|\nabla_t f(x_*) - \nabla_t f(x_n^j)\|^2 &\leq& f(x_*)-f(x_n^j)-\ve{\nabla f(x_n^j)}{x_*-x_n^j},
\end{eqnarray*}
for $i,j,t=1,2,...,p$ and $m,n=0,1,...,N$, that is
\begin{eqnarray*}
\frac{1}{2L_t}\|\mathbf{U}_t^\top\nabla f(x_m^i) - \mathbf{U}_t^\top\nabla f(x_n^j)\|^2 &\leq& f(x_m^i)-f(x_n^j)-\ve{\nabla f(x_n^j)}{x_m^i-x_n^j},\\
\frac{1}{2L_t}\|\mathbf{U}_t^\top\nabla f(x_m^i)\|^2 &\leq& f(x_m^i)-f(x_*),\\
\frac{1}{2L_t}\| \mathbf{U}_t^\top\nabla f(x_n^j)\|^2 &\leq& f(x_*)-f(x_n^j)-\ve{\nabla f(x_n^j)}{x_*-x_n^j},
\end{eqnarray*}
where we use the fact $\mathbf{U}_t^\top\nabla f(x)=\nabla_t f(x)$.

In this paper we deal with a standard case to get the insights, here we assume all the block partial Lipschitz constants are equal, that is $L_1=L_2=...=L_p=L_c$.
We define
\begin{align*}
    \delta_{k,i} &:= \frac{1}{pL_cR^2(x_0)}(f(x_k^i)-f(x_\ast)), \\
    g_{k,i} &:= \frac{1}{pL_cR(x_0)}\nabla f(x_k^i), \\
    \delta_{*} &:= \frac{1}{pL_cR^2(x_0)}(f(x_\ast)-f(x_\ast))=0, \\
    g_{*} &:= \frac{1}{pL_cR(x_0)}\nabla f(x_\ast)=0,
\end{align*}
for every $k=0,\ldots, N,~i=1,\ldots, p$.
In view of Algorithm~(\ref{alg:bcd}), since $x_k^p=x_{k+1}^0$ for $k=0,...,N$, obviously we have
\begin{align}
\delta_{k,p}=\delta_{k+1,0}\quad \text{and} \quad g_{k,p}=g_{k+1,0}.
\end{align}

In view of the above notations, Problem \eqref{P:baseproblem} can now be relaxed by discarding the constrains $f(x)\in \mathcal{F}_{\mathbf{L},\mathbf{U}}(\mathbb{R}^D)$ to the following form
\begin{align*}
\begin{aligned}
    \max_{\substack{x_k^i\in\mathbb{R}^D,~g_{k,i} \in\mathbb{R}^D,~\delta_{k,i} \in\mathbb{R}, \\ k=0,\ldots, N,~i=1,\ldots, p.}}\ & pL_cR^2(x_0) \delta_{N,p} \\
    \text{s.t. }\quad
                        & x_k^i=x_k^{i-1}-p R(x_0)\mathbf{U}_i \mathbf{U}_i^\top g_{k,i-1},,\\
                        & \delta_{k,p}=\delta_{k+1,0},\quad g_{k,p}=g_{k+1,0}, \\
                        & \frac{p}{2}\| \mathbf{U}_t^\top g_{m,i} - \mathbf{U}_t^\top g_{n,j} \| ^2 \leq \delta_{m,i} - \delta_{n,j} - \frac{ \langle g_{n,j},  x_m^i-x_n^j  \rangle}{R(x_0)}, \\
                        & \frac{p}{2}\| \mathbf{U}_t^\top g_{m,i}  \| ^2 \leq \delta_{m,i}, \\
                        & \frac{p}{2}\| \mathbf{U}_t^\top g_{n,j} \| ^2 \leq - \delta_{n,j} - \frac{ \langle g_{n,j},  x_*-x_n^j  \rangle}{R(x_0)}, \\
                        & \frac{p}{2}\| \mathbf{U}_t^\top g_{0,0} \| ^2 \leq - \delta_{0,0} - \frac{ \langle g_{0,0},  x_*-x_0^0  \rangle}{R(x_0)}, \\
                        & \quad  m,n,k=0,\ldots, N,\quad i,j,t=1,\ldots, p.
     \end{aligned}
        \tag{P1}\label{P:baseproblem_equal}
\end{align*}

We try to relax the above problem, if we set $m=n=k$ and $j=i+1$, then we have
\begin{align*}
\begin{aligned}
    \max_{\substack{x_k^i\in\mathbb{R}^D,~g_{k,i} \in\mathbb{R}^D,~\delta_{k,i} \in\mathbb{R}, \\ k=0,\ldots, N,~i=1,\ldots, p.}}\ & pL_cR^2(x_0) \delta_{N,p} \\
    \text{s.t. }\quad
                        & x_k^i=x_k^{i-1}-p R(x_0)\mathbf{U}_i \mathbf{U}_i^\top g_{k,i-1},\\
                        & \delta_{k,p}=\delta_{k+1,0},\quad g_{k,p}=g_{k+1,0}, \\
                        & \frac{p}{2}\| \mathbf{U}_t^\top g_{k,i-1} - \mathbf{U}_t^\top g_{k,i} \| ^2 \leq \delta_{k,i-1} - \delta_{k,i} - \frac{ \langle g_{k,i},  x_k^{i-1}-x_k^i  \rangle}{R(x_0)}, \\
                        & \frac{p}{2}\| \mathbf{U}_t^\top g_{k,i}  \| ^2 \leq \delta_{k,i}, \\
                        & \frac{p}{2}\| \mathbf{U}_t^\top g_{k,i} \| ^2 \leq - \delta_{k,i} - \frac{ \langle g_{k,i},  x_*-x_k^i  \rangle}{R(x_0)}, \\
                        & \frac{p}{2}\| \mathbf{U}_t^\top g_{0,0} \| ^2 \leq - \delta_{0,0} - \frac{ \langle g_{0,0},  x_*-x_0^0  \rangle}{R(x_0)}, \\
                        & \quad k=0,\ldots, N,\quad i,t=1,\ldots, p.
     \end{aligned}
        \tag{P2}\label{P:prime_problem2}
\end{align*}

Same as in~\cite{drori2012performance}, the Problem~\eqref{P:prime_problem2} is invariant under the transformation $g_{k,i}^\prime
\leftarrow Q g_{k,i}$, $x_{k,i}^\prime \leftarrow Q x_{k,i}$ for any orthogonal transformation $Q$. We can
therefore assume without loss of generality that $x_\ast- x_0= \|x_\ast-x_0\|\nu$, where $\nu$
is any given unit vector in $\mathbb{R}^D$. Therefore, we have
\begin{align*}
\frac{p}{2}\| \mathbf{U}_t^\top g_{k,i} \| ^2 \leq - \delta_{k,i} - \frac{ \langle g_{k,i},  \|x_*-x_0\|\nu+x_0-x_{k,i} \rangle}{R(x_0)}.
\end{align*}
and
\begin{align*}
x_0-x_{k,i} = pR(x_0) \sum_{\{k',i':k'p+i'\leq kp+i\}}\mathbf{U}_{i'} \mathbf{U}_{i'}^\top g_{k',i'-1}
\end{align*}

In order to simplify notation, we denote $\sum_{\{k',i':k'p+i'\leq kp+i\}}$ as $\sum_{k'p+i'\leq kp+i}$ in the following
Now we can remove some constraints from Problem \eqref{P:prime_problem2} to further simplify the analysis:
\begin{align*}
\begin{aligned}
    \max_{\substack{x_k^i\in\mathbb{R}^D,~g_{k,i} \in\mathbb{R}^D,~\delta_{k,i} \in\mathbb{R}, \\ k=0,\ldots, N,~i=1,\ldots, p.}}\ & pL_cR^2(x_0) \delta_{N,p} \\
    \text{s.t. }\quad
                        & \frac{p}{2}\| \mathbf{U}_t^\top g_{k,i-1} - \mathbf{U}_t^\top g_{k,i} \| ^2 \leq \delta_{k,i-1} - \delta_{k,i} - p \langle g_{k,i}, \mathbf{U}_i \mathbf{U}_i^\top g_{k,i-1} \rangle, \\
                        & \frac{p}{2}\| \mathbf{U}_t^\top g_{k,i} \| ^2 \leq - \delta_{k,i} - \langle g_{k,i},  \alpha \nu + p \sum_{k'p+i'\leq kp+i}\mathbf{U}_{i'} \mathbf{U}_{i'}^\top g_{k',i'-1}   \rangle, \\
                        & \frac{p}{2}\| \mathbf{U}_t^\top g_{0,0} \| ^2 \leq - \delta_{0,0} - \langle g_{0,0},  \alpha \nu \rangle, \\
                        & \delta_{k,p}=\delta_{k+1,0},\quad g_{k,p}=g_{k+1,0}, \\
                        & k=0,\ldots, N,\quad i,t=1,\ldots, p,
     \end{aligned}
        \tag{P3}\label{P:prime_problem3}
\end{align*}
where $\alpha = \frac{\|x_*-x_0\|}{R(x_0)}$. It is obvious that $\alpha \leq 1$ from the definition~\eqref{eq:rx0}.

Let $G$ denote the $((N+1)p+1) \times D$ matrix whose rows are $g_{0,0}^\top, g_{0,1}^\top, ...,g_{k,i}^\top,...,g_{N,p}^\top $, and $u_{k,i}\in \mathbb{R}^{(N+1)p+1}$ be $e_{kp+i+1}$, the ($kp+i+1$)th standard unit vector for $i=1,2,...,p$ and $k=0,...,N$,. Then we have
\[
   \mathbf{U}_t^\top g_{k,i}= \mathbf{U}_t^\top G^\top u_{k,i},\; \trace (\mathbf{U}_t^\top G^\top u_{m,i} u_{n,j}^\top G \mathbf{U}_t)=\langle \mathbf{U}_t^\top g_{m,i}, \mathbf{U}_t^\top g_{n,j} \rangle, \;\mbox{and}\; \langle G^\top u_{k,i},\nu \rangle=\langle g_{k,i}, \nu\rangle
\]
for any $m,n,k,i,j,t$. Let $\delta = (\delta_{0,0},\delta_{0,1},...,\delta_{k,i},...,\delta_{N,p} ) \in \mathbb{R}^{(N+1)p+1}$

Let $t=i$ Problem \eqref{P:prime_problem3}, then it can be transformed into a more compact form in terms of $G$ and $\delta$
\begin{align*}
\begin{aligned}
    \min_{\substack{G \in \mathbb{R}^{((N+1)p+1) \times D}, \\ \delta  \in \mathbb{R}^{(N+1)p+1}}}\ & - \delta_{N,p} \\
    \text{s.t. }\quad
                        & \frac{p}{2}\trace (\mathbf{U}_i\mathbf{U}_i^\top G^\top (u_{k,i-1} u_{k,i-1}^\top + u_{k,i} u_{k,i}^\top) G )\leq \delta_{k,i-1} - \delta_{k,i} \\
                        & \frac{p}{2}\trace (\mathbf{U}_i\mathbf{U}_i^\top G^\top  u_{k,i} u_{k,i}^\top G )
                         + \trace ( \alpha \nu u_{k,i}^\top G ) \\ & \quad \quad + \frac{p}{2} \sum_{k'p+i'\leq kp+i} \trace (\mathbf{U}_{i'}\mathbf{U}_{i'}^\top G^\top ( u_{k,i} u_{k',i'-1}^\top+u_{k',i'-1}u_{k,i}^\top) G )  \leq - \delta_{k,i},  \\
                        & \frac{p}{2}\trace (\mathbf{U}_i\mathbf{U}_i^\top G^\top  u_{0,0} u_{0,0}^\top G ) + \trace ( \alpha \nu u_{0,0}^\top G )  \leq - \delta_{0,0},  \\
                        & \delta_{k,p}=\delta_{k+1,0},\quad g_{k,p}=g_{k+1,0}, \\
                        & k=0,\ldots, N,\quad i=1,\ldots, p,
     \end{aligned}
        \tag{P4}\label{P:prime_G}
\end{align*}
where in order for convenience,  we recast the above as a minimization problem, and we also omit the
fixed term $pL_cR^2(x_0)$ from the objective.

Attaching the dual multipliers $$\lambda :=(\lambda_{0,1},...,\lambda_{k,i},...,\lambda_{N,p})^\top \in \mathbb{R}^{(N+1)p}_{+}$$ and
$$\tau:=(\tau_{0,0},...,\tau_{k,i},...,\tau_{N,p})^\top \in \mathbb{R}^{(N+1)p+1}_{+}$$ to the first and second set of
inequalities respectively, and using the notation $$\delta=(\delta_{0,0},...,\delta_{k,i},...,\delta_{N,p})$$, we
get that the Lagrangian of this problem is given as a sum of two separable functions in the
variables $(\delta, G)$:
\begin{eqnarray*}
    L(G,\delta, \lambda, \tau)&=&
         -\delta_{N,p}+ \sum_{k=0}^{N}\sum_{i=1}^{p} \lambda_{k,i}(\delta_{k,i} - \delta_{k,i-1}) + \sum_{k=0}^{N}\sum_{i=1}^{p} \tau_{k,i} \delta_{k,i} + \tau_{0,0}\delta_{0,0}   \\
         & &+ \frac{p}{2} \sum_{k=0}^{N}\sum_{i=1}^{p}  \lambda_{k,i} \trace (\mathbf{U}_i\mathbf{U}_i^\top G^\top (u_{k,i-1} u_{k,i-1}^\top + u_{k,i} u_{k,i}^\top) G )  \\
         & &+ \sum_{k=0}^{N}\sum_{i=1}^{p} \tau_{k,i} [\frac{p}{2}\trace (\mathbf{U}_i\mathbf{U}_i^\top G^\top u_{k,i}u_{k,i}^\top G)+\trace (\alpha \nu  u_{k,i}^\top G) \\
         & &+ \frac{p}{2} \sum_{k'p+i'\leq kp+i}  \trace  ( \mathbf{U}_{i'} \mathbf{U}_{i'}^\top G^\top (u_{k',i'-1}u_{k,i}^\top + u_{k,i}u_{k',i'-1}^\top) G)]\\
        & &+  \tau_{0,0} [\frac{p}{2}\trace (\mathbf{U}_i\mathbf{U}_i^\top G^\top u_{0,0}u_{0,0}^\top G)+\trace (\alpha \nu  u_{0,0}^\top G) ]\\
        & \equiv & L_1(\delta, \lambda, \tau)+ L_2(G,\lambda, \tau).
\end{eqnarray*}

The dual objective function is then defined by $$H(\lambda, \tau)=\min_{G, \delta} L(G,\delta,
\lambda\, \tau)=\min_\delta L_1(\delta, \lambda, \tau)+ \min_{G} L_2(G, \lambda, \tau),$$
and the dual problem of Problem \eqref{P:prime_G} is then given by
\begin{align*}
    \begin{aligned}
        \max \{ H(\lambda, \tau): \lambda \in \mathbb{R}^{(N+1)p+1}_{+}, \tau \in \mathbb{R}^{(N+1)p+1}_{+}\}.
    \end{aligned}
\end{align*}
Since $L_1(\cdot,\lambda, \tau)$ is linear in $\delta$, we have $\min_\delta
L_1(\delta,\lambda, \tau) =0$ whenever
\begin{eqnarray}
            -\lambda_{0,1} +\tau_{0,0} &= & 0, \nonumber\\
            \lambda_{k,i}-\lambda_{k,i+1}+\tau_{k,i} &=& 0, \quad (k=1,\dots,N,\quad i=1,...,p-1), \label{eqtau}\\
             -1+\lambda_{N,p}+\tau_{N,p} &=& 0, \nonumber
\end{eqnarray}
and $-\infty$ otherwise.

According to Lemma~\ref{lem:quadratic_min}, we have
\begin{align*}
     \min_{G \in \mathbb{R}^{((N+1)p+1) \times D } } L_2(G,\lambda, \tau) =
     \min_{w \in \mathbb{R}^{(N+1)p+1 } } L_2(w \nu^\top,\lambda, \tau).
\end{align*}

Let $\nu$ be $(\frac{1}{\sqrt{D}},...,\frac{1}{\sqrt{D}})^\top$, then we have $\nu^\top\mathbf{U}_i\mathbf{U}_i^\top  \nu = \frac{D_i}{D}$.
Therefore for any $(\lambda,\tau)$ satisfying~\eqref{eqtau}, we have obtained that the dual objective is upper bounded by
\begin{eqnarray*}
  H(\lambda, \tau)&\leq & =  \min_{w \in \mathbb{R}^{(N+1)p+1 } } L_2(w \nu^\top,\lambda, \tau) \\
  &=& \min_{w \in \mathbb{R}^{(N+1)p+1 } }\{ \frac{p}{2} \sum_{k=0}^{N}\sum_{i=1}^{p} \frac{D_i}{D}\lambda_{k,i}  w^\top (u_{k,i-1} u_{k,i-1}^\top + u_{k,i} u_{k,i}^\top)w  \\
  &+&  \sum_{k=0}^{N}\sum_{i=1}^{p} \tau_{k,i} [\frac{p}{2}\frac{D_i}{D} w^\top u_{k,i}u_{k,i}^\top w+ \alpha   u_{k,i}^\top w + \frac{p}{2}\sum_{k'p+i'\leq kp+i}  \frac{D_{i'}}{D} w^\top( u_{k,i}u_{k',i'-1}^\top+ u_{k',i'-1}u_{k,i}^\top) w] \\
  &+& \tau_{0,0} [\frac{p}{2}\frac{D_i}{D} w^\top u_{0,0}u_{0,0}^\top w+\alpha   u_{0,0}^\top w ] \}\\
        &=&  \max_{t\in \mathbb{R}} \{-\frac{1}{2}t :
        w^\top A w + \alpha \tau^\top w \geq -\frac{1}{2}t,\ \forall w\in \mathbb{R}^{(N+1)p+1 }
        \}\\
 &=& \max_{t\in \mathbb{R}} \left\{
            -\frac{1}{2}t:
            \begin{pmatrix}
                  A  & \frac{1}{2} \tau \\
                \frac{1}{2}\tau^\top       & \frac{1}{2}t
        \end{pmatrix} \succeq 0
        \right\},
\end{eqnarray*}
where
\begin{eqnarray*}
A&=&\frac{p}{2} \sum_{k=0}^{N}\sum_{i=1}^{p} \frac{D_i}{D} \lambda_{k,i}  (u_{k,i-1} u_{k,i-1}^\top + u_{k,i} u_{k,i}^\top)
  \\ &+&  \sum_{k=0}^{N}\sum_{i=1}^{p} \tau_{k,i} [\frac{p}{2}\frac{D_i}{D} u_{k,i}u_{k,i}^\top + \frac{p}{2}\sum_{k'p+i'\leq kp+i}  \frac{D_{i'}}{D}  (u_{k,i}u_{k',i'-1}^\top+u_{k',i'-1}u_{k,i}^\top)]
  + \tau_{0,0} \frac{p}{2}\frac{D_i}{D} u_{0,0}u_{0,0}^\top.
\end{eqnarray*}

If all the block have equal size, that is $\frac{D_i}{D}=\frac{1}{p}$ for every $i=1,...,p$, then we get
\begin{eqnarray*}
A&=&\frac{1}{2} \sum_{k=0}^{N}\sum_{i=1}^{p} \lambda_{k,i}  (u_{k,i-1} u_{k,i-1}^\top + u_{k,i} u_{k,i}^\top) \\
  &+&  \sum_{k=0}^{N}\sum_{i=1}^{p} \tau_{k,i} [\frac{1}{2} u_{k,i}u_{k,i}^\top + \sum_{k'p+i'\leq kp+i} \frac{1}{2} (u_{k,i}u_{k',i'-1}^\top+u_{k',i'-1}u_{k,i}^\top)]
  + \tau_{0,0} \frac{1}{2} u_{0,0}u_{0,0}^\top. \label{matrix_A}
\end{eqnarray*}

Now we obtain an upper bound for the optimal value of Problem~\eqref{P:prime_problem3}:
\begin{align*}
\begin{aligned}
 \min_{t\in \mathbb{R}} \left\{
            \frac{1}{2}pL_cR^2(x_0) t :
            \begin{pmatrix}
                  A  & \frac{1}{2} \tau \\
                \frac{1}{2}\tau^\top       & \frac{1}{2}t
        \end{pmatrix} \succeq 0
        \right\}.
 \end{aligned}
        \tag{D}\label{D:dual_problem1}
\end{align*}

\section{New bound of BCD}

Note~\eqref{eqtau} and~\eqref{matrix_A}, we have
\begin{eqnarray*}
\tau &=& (\tau_{0,0},\tau_{0,1},...,\tau_{k,i},...,\tau_{N,p})^\top \\
&=& (\lambda_{0,0},\lambda_{0,2}-\lambda_{0,1},...,\lambda_{k,i+1}-\lambda_{k,i},...,1-\lambda_{N,p})^\top
\end{eqnarray*}
and $2A$ becomes
\begin{equation}
\left(
  \begin{array}{ccccccc}
   2\lambda_{0,1}& \lambda_{0,2}-\lambda_{0,1} & \cdots & \lambda_{k,i+1}-\lambda_{k,i}  & \cdots & \lambda_{N,p}-\lambda_{N,p-1} & 1- \lambda_{N,p}\\
    \lambda_{0,2}-\lambda_{0,1} & 2\lambda_{0,2} &  & \lambda_{k,i+1}-\lambda_{k,i}  &  & \lambda_{N,p}-\lambda_{N,p-1} & 1- \lambda_{N,p}\\
   \vdots &  & \ddots & & &  & \vdots\\
    \lambda_{k,i+1}-\lambda_{k,i} & \lambda_{k,i+1}-\lambda_{k,i} &  &2\lambda_{k,i} &  & \lambda_{N,p}-\lambda_{N,p-1} & 1- \lambda_{N,p}\\
    \vdots & &  &  & \ddots & & \vdots\\
    \lambda_{N,p}-\lambda_{N,p-1} & \lambda_{N,p}-\lambda_{N,p-1} & & \lambda_{N,p}-\lambda_{N,p-1}& & 2\lambda_{N,p} & 1- \lambda_{N,p}\\
    1- \lambda_{N,p}&   1- \lambda_{N,p} &   \cdots &   1- \lambda_{N,p} &  \cdots & 1- \lambda_{N,p} & 1
  \end{array}
\right).
\end{equation}

According to Appendix~\ref{sec::app3}, if we set
\begin{eqnarray*}
\lambda_{k,i} &=& \frac{kp+i}{2(N+1)p+1-kp-i}, \quad k=0,...,N, \quad p=1,...,p, \\
t &=& \frac{1}{2(N+1)p+1}.
\end{eqnarray*}
we have $\begin{pmatrix}
                  2A  &  \tau \\
               \tau^\top       & t
        \end{pmatrix} \succeq 0.$
Thus we have the following new upper bound on the complexity of the BCD:
\begin{theorem}\label{T:gradUpperBound}
Let $f(x)\in  \mathcal{F}_{\mathbf{L},\mathbf{U}}(\mathbb{R}^D)$ and let $x_0,\dots,x_N\in \mathbb{R}^D$ be generated by
Algorithm~\ref{alg:bcd} with $L_1=L_2=...=L_p=L_c$ and $D_1=D_2=...=D_p$. Then we have
\begin{equation}\label{E:gmbound}
    f(x_N)-f(x_*)\leq \frac{1}{4(N+1)p+2}pL_cR^2(x_0).
\end{equation}
\end{theorem}

\begin{remark}
From  above theorem, we notice that our bound is $16p^3$ times smaller than for the known bound~\eqref{eq:bcd_conv_beck} (with $L_1=L_2=...=L_p=L_c$)
\begin{align*}
f(x_k)-f(x_*)\leq 4L_c(1+p^3)R^2(x_0)\frac{1}{k+8/p}.
\end{align*}
\end{remark}

\section{Numerical test}
\label{sec:num_test}

Consider the least squares problem
\begin{align}
  \min_{x \in \mathbb{R}^D} \,\frac{1}{2}\|\textbf{A}x-b\|^2,
  \label{eq:problem_least_square}
\end{align}
where $\textbf{A} \in \mathbb{R}^{D\times D}$, $b \in \mathbb{R}^D$. $\textbf{A}$ is a nonsingular matrix, so obviously the optimal solution of the problem is the vector $\textbf{A}^{-1}b$ and the optimal value is $f_*=0$. We consider the partition of the variables to $p$ blocks, each with $n/p$ variables (we assume that $p$ divides $n$). We will also use the notation
$$\mathbf{A}=(\textbf{A}_1\quad \textbf{A}_2 \dots \quad \textbf{A}_p)$$
where $\textbf{A}_i$ is the submatrix of $\textbf{A}$ comprising the columns corresponding to the $i$-th block,
that is, columns $(i-1)n/p+1,(i-1)n/p+2,...,in/p$.

We consider $n=100$, and three choices of $p$: 2,5, 20, and 100. The results together with classical bound on the convergence rate of the sequence of the BCD method are summarized in Figure 1.

\section{Conclusion}
This paper provide a novel and better analytical convergence bound, that is  $16p^3$ times as small as previous best, for the sequence of BCD methods for unconstrained smooth convex functions. Extending this approach to general, such randomized BCD type method or stochastic gradient method, is important future work. In a broader context, we believe that the current paper could serve as a basis for examining the method on the PEP approach to various BCD related methods.

% ---- Bibliography ----
%\bibliographystyle{abbrv}  %this one
%\bibliography{IEEEabrv,CommunityBIB-Jerry.bib}
\bibliography{bcd_pep_arxiv}

\appendix

\section{Proof of Lemma~\ref{lem:bound_diff_nabla}}

In this appendix, we complete the proofs of the Lemma~\ref{lem:bound_diff_nabla}.

For all  $h^i \in \mathbb{R}^{D_i}$ and $x \in \mathbb{R}^{D}$, we have
\begin{eqnarray*}
  f(x+\mathbf{U}_i h_i)
   &=&  f(x) + \int_0^1 \ve{\nabla f(x+\theta \mathbf{U}_i h_i)}{\mathbf{U}_i h_i} d\theta\\
  &=& f(x) + \int_0^1 \ve{\mathbf{U}_i^\top \nabla f(x+\theta \mathbf{U}_i h_i)}{h_i} d\theta \\
   &=&  f(x) +\ve{ \nabla_i f(x)}{h_i}+ \int_0^1 \ve{ \nabla_i f(x+\theta \mathbf{U}_i h_i)-\nabla_i f(x)}{h_i}  d\theta \\
    &\leq&  f(x) +\ve{ \nabla_i f(x)}{h_i}+ \int_0^1 \| \ve{ \nabla_i f(x+\theta \mathbf{U}_i h_i)-\nabla_i f(x)}{h_i} \| d\theta \\
    &\leq&  f(x) +\ve{ \nabla_i f(x)}{h_i}+ \int_0^1 \| { \nabla_i f(x+\theta \mathbf{U}_i h_i)-\nabla_i f(x)} \| \|{h_i} \| d\theta \\
      &\leq&  f(x) +\ve{ \nabla_i f(x)}{h_i}+ \int_0^1 \theta L_i\| h_i\|^2 d\theta  \\
      &=&  f(x) +\ve{ \nabla_i f(x)}{h_i}+ \frac{L_i}{2}\| h_i\|^2,
\end{eqnarray*}
where the second inequality follows from the Cauchy-Schwartz inequality and the third inequality follows from~\eqref{eq:block_coordinate_lipschitz}.
In short from above we have
\begin{equation}\label{eq:bound_gradient_descent}
 f(x+\mathbf{U}_i h_i) \leq f(x) +\ve{ \nabla_i f(x)}{h_i}+ \frac{L_i}{2}\| h_i\|^2.
\end{equation}

Then consider the function $\phi(y)=f(y)-\ve{\nabla f(x_0)}{y}$. The gradient of $\phi(y)$ is $\nabla f(y) - \nabla f(x_0)$, which is obvious block-coordinate-wise Lipschitz continuous with constants $L_i$, thais belong to the class $\mathcal{F}_{\mathbf{L},\mathbf{U}}(\mathbb{R}^D)$, same as $f(x)$, and $x_0$ is one of its optimal points. Therefor in view of~\eqref{eq:bound_gradient_descent}, we get
\begin{equation}\label{eq:bound_diff_nabla_phi_y}
\phi(y+\mathbf{U}_i h_i) \leq \phi(y) +\ve{ \nabla_i \phi(y)}{h_i}+ \frac{L_i}{2}\| h_i\|^2.
\end{equation}

Let $h_i=-\frac{1}{L_i}\nabla_i \phi(y)$ in~\eqref{eq:bound_diff_nabla_phi_y}, we have
\begin{eqnarray*}
  \phi(x_0)
   &\leq&  \phi(y+\mathbf{U}_i(-\frac{1}{L_i}\nabla_i \phi(y)))\\
      &\leq&  \phi(y) +\ve{ \nabla_i \phi(y)}{-\frac{1}{L_i}\nabla_i \phi(y)}+ \frac{L_i}{2}\| -\frac{1}{L_i}\nabla_i \phi(y)\|^2 \\
      &=& \phi(y) - \frac{1}{2L_i} \|\nabla_i \phi(y)\|^2,
\end{eqnarray*}
where the first inequality follows from the fact that $x_0$ is one of $\phi(y)$'s optimal points.

Since $\phi(x_0)=f(x_0)-\ve{\nabla f(x_0)}{x_0}$ and $\nabla_i \phi(y)=\nabla_i f(y) - \nabla_i f(x_0)$, thus we have
\begin{equation}
f(x_0)-\ve{\nabla f(x_0)}{x_0} \leq f(y)-\ve{\nabla f(x_0)}{y} - \frac{1}{2L_i}\|\nabla_i f(y) - \nabla_i f(x_0)\|^2. \nonumber
\end{equation}

Therefore we get (\ref{eq:bound_diff_nabla_i}) and the lemma is proved.

\section{Proof of Lemma~\ref{lem:quadratic_min}}

First, for any $\xi \in \mathbb{R}^n$, we have $f(\xi b^\top)=b^\top A b\xi^\top B\xi + 2\|b\|^2a^\top \xi$. Thus, likewise,
 $\inf\{f(\xi b^\top): \xi \in \mathbb{R}^n\} >-\infty$ if and only if both $A \succeq 0$ and $B \succeq 0$, and there exists ${\bar \xi}\in \mathbb{R}^n$ such
that
\begin{equation}
\label{s2} b^\top A b B{\bar \xi} + \|b\|^2a=0,
\end{equation}
and using (\ref{s2}) it
follows $\inf_\xi f(\xi b^\top)=f({\bar \xi}b^\top)=\|b\|^2a^\top{\bar \xi}=\trace (ba^\top{\bar
\xi}b^\top)$.

Now, we recall that $\inf \{f(X): X \in \mathbb{R}^{n\times m}\}
>-\infty$ if and only if both $A \succeq 0$ and $B \succeq 0$, and there exists at least one solution ${\bar X}$ such
that
\begin{equation}\label{s1}
 A{\bar X}^\top B +ba^\top=0,
\end{equation}
i.e., the above is just $\nabla f (X)= 0$ and characterizes the minimizers of the convex
function $f(X)$. Using (\ref{s1}) it follows that $\inf_X f(X)=f({\bar X})=\trace (ba^\top{\bar
X})$.

Now, using (\ref{s1})-(\ref{s2}), one obtains $A{\bar X}^\top B=-ba^\top$ and $A({\bar X}-
{\bar \xi}b^\top)B=0$, and hence it follows that
\begin{eqnarray*}
f({\bar X})-f({\bar \xi})&=&\trace (ba^\top({\bar X}- {\bar \xi}b^\top))\\
                        &=&\trace (-A{\bar X}^\top B({\bar X}- {\bar \xi}b^\top)=0.
                        \end{eqnarray*}

\section{Positive of $A$}
\label{sec::app3}

Recall that $2A$ equals
\begin{equation}
\left(
  \begin{array}{ccccccc}
   2\lambda_{0,1}& \lambda_{0,2}-\lambda_{0,1} & \cdots & \lambda_{k,i+1}-\lambda_{k,i}  & \cdots & \lambda_{N,p}-\lambda_{N,p-1} & 1- \lambda_{N,p}\\
    \lambda_{0,2}-\lambda_{0,1} & 2\lambda_{0,2} &  & \lambda_{k,i+1}-\lambda_{k,i}  &  & \lambda_{N,p}-\lambda_{N,p-1} & 1- \lambda_{N,p}\\
   \vdots &  & \ddots & & &  & \vdots\\
    \lambda_{k,i+1}-\lambda_{k,i} & \lambda_{k,i+1}-\lambda_{k,i} &  &2\lambda_{k,i} &  & \lambda_{N,p}-\lambda_{N,p-1} & 1- \lambda_{N,p}\\
    \vdots & &  &  & \ddots & & \vdots\\
    \lambda_{N,p}-\lambda_{N,p-1} & \lambda_{N,p}-\lambda_{N,p-1} & & \lambda_{N,p}-\lambda_{N,p-1}& & 2\lambda_{N,p} & 1- \lambda_{N,p}\\
    1- \lambda_{N,p}&   1- \lambda_{N,p} &   \cdots &   1- \lambda_{N,p} &  \cdots & 1- \lambda_{N,p} & 1
  \end{array}
\right) \nonumber
\end{equation}
for
\begin{eqnarray*}
\lambda_{k,i} &=& \frac{kp+i}{2(N+1)p+1-kp-i}, \quad k=0,...,N, \quad p=1,...,p.
\end{eqnarray*}

We begin by deriving a recursion rule for the determinant of matrices of the following form:
\[
    A_{k,i}=
    \begin{pmatrix}
        2\lambda_{0,1}& \lambda_{0,2}-\lambda_{0,1} & \cdots & \lambda_{k,i+1}-\lambda_{k,i}   \\
         \lambda_{0,2}-\lambda_{0,1} & 2\lambda_{0,2} &  & \lambda_{k,i+1}-\lambda_{k,i} \\
         \vdots &  & \ddots &  \\
        \lambda_{k,i+1}-\lambda_{k,i} & \lambda_{k,i+1}-\lambda_{k,i} &  &2\lambda_{k,i}
    \end{pmatrix}.
\]

To find the determinant of $A_{k,i}$,
subtract the one before last row multiplied by $\frac{\lambda_{k,i+1}-\lambda_{k,i}}{\lambda_{k,i}-\lambda_{k,i-1}}$ from the last row: the last row becomes
\[
    (0,\dots,0,\lambda_{k,i+1}-\lambda_{k,i}-\frac{\lambda_{k,i+1}-\lambda_{k,i}}{\lambda_{k,i}-\lambda_{k,i-1}}2\lambda_{k,i-1},2\lambda_{k,i}-\frac{\lambda_{k,i+1}-\lambda_{k,i}}{\lambda_{k,i}-\lambda_{k,i-1}}(\lambda_{k,i+1}-\lambda_{k,i})).
\]
Expanding the determinant along the last row we get
\begin{align*}
    \det A_{k,i} =& (2\lambda_{k,i}-\frac{\lambda_{k,i+1}-\lambda_{k,i}}{\lambda_{k,i}-\lambda_{k,i-1}}(\lambda_{k,i+1}-\lambda_{k,i})) \det A_{k,i-1}\\
    -&(\lambda_{k,i+1}-\lambda_{k,i}-\frac{\lambda_{k,i+1}-\lambda_{k,i}}{\lambda_{k,i}-\lambda_{k,i-1}}2\lambda_{k,i-1}) \det (A_{k,i-1})_{kp+i,kp+i-1}
\end{align*}
where $  (A_{k,i-1})_{kp+i,kp+i-1}$ denotes the $kp+i,kp+i-1$ minor:
\[
    (A_{k,i-1})_{kp+i,kp+i-1}=
    \begin{pmatrix}
        2\lambda_{0,1}& \lambda_{0,2}-\lambda_{0,1} & \cdots & \lambda_{k,i+1}-\lambda_{k,i}   \\
         \lambda_{0,2}-\lambda_{0,1} & 2\lambda_{0,2} &  & \lambda_{k,i+1}-\lambda_{k,i} \\
         \vdots &  & \ddots &  \\
        \lambda_{k,i}-\lambda_{k,i-1} & \lambda_{k,i}-\lambda_{k,i-1} &  & \lambda_{k,i+1}-\lambda_{k,i}
    \end{pmatrix}.
    .
\]
If we multiply the last column of $(A_{k,i-1})_{kp+i,kp+i-1}$ by $\frac{\lambda_{k,i}-\lambda_{k,i-1}}{\lambda_{k,i+1}-\lambda_{k,i}}$
we get a matrix that is different from $A_{k,i-1}$ by only the corner element.
Thus by basic determinant properties we get that
\[
    \frac{\lambda_{k,i}-\lambda_{k,i-1}}{\lambda_{k,i+1}-\lambda_{k,i}}\det (A_{k,i-1})_{kp+i,kp+i-1} = \det A_{k,i-1}+( \lambda_{k,i+1}-\lambda_{k,i}-2\lambda_{k,i-1}) \det A_{k,i-2}.
\]
Combining these two results, we have found the following recursion rule for $\det A_{k,i}$, $kp+i\geq 2$:
\begin{align}\label{E:mrecursion}
    \det A_{k,i}=&\left(2\lambda_{k,i}-\frac{2 ( \lambda_{k,i}-\lambda_{k,i-1})^2}{ \lambda_{k,i-1}-\lambda_{k,i-2}}+\frac{( \lambda_{k,i}-\lambda_{k,i-1})^2 2\lambda_{k,i-1}}{( \lambda_{k,i-1}-\lambda_{k,i-2})^2} \right) \det A_{k,i-1}
    \nonumber\\ -&( \lambda_{k,i}-\lambda_{k,i-1})^2\left(1-\frac{2\lambda_{k,i-1}}{ \lambda_{k,i-1}-\lambda_{k,i-2}}\right)^2 \det A_{k,i-2} .
\end{align}
Obviously, the recursion base cases are given by
\begin{align*}
    & \det A_{0,1} =  2\lambda_{0,1}, \\
    & \det A_{0,2} =  4\lambda_{0,1}\lambda_{0,2} - (\lambda_{0,2}-\lambda_{0,1})^2.
\end{align*}

Since
\begin{eqnarray*}
\lambda_{k,i} &=& \frac{kp+i}{2(N+1)p+1-kp-i}, \quad k=0,...,N, \quad p=1,...,p.
\end{eqnarray*}
we get that $A_{k,i}$ is the $kp+i+1$'th leading principal minor of the matrix $A$.
The recursion rule~\eqref{E:mrecursion} can now be solved. The solution is given by:
\begin{align}
    \det A_{k,i} =& \frac{(2(N+1)p+1)^2}{(2(N+1)p-kp-i)^2} ( 1 \nonumber  \\
    +&\sum_{k'p+i'\leq kp+i-1} \frac{2(N+1)p-2k'p-2i'-1}{2(N+1)p+4(N+1)p (k'p+i') - 2 (k'p+i')^2+1}) \label{E:mkdet} \\
    \times¡¡&\prod_{k'p+i'\leq kp+i-1} \frac{2(N+1)p+4(N+1)p(k'p+i')-2(k'p+i')^2+1}{(2(N+1)p+1-(k'p+i'))^2},  \nonumber
\end{align}
for $k=0,\dots,N-1$ and $i=1,...,p$, and
\begin{align}
    \det A_{N,p}
        =& \frac{(2N+1)^2}{(N+1)^2}\prod_{i=0}^{N-1} \frac{2N+4Ni-2i^2+1}{(2N+1-i)^2} \label{E:mndet}.
\end{align}

We now proceed to verify the expressions \eqref{E:mkdet} and \eqref{E:mndet} given above. We will show that these expressions satisfy the recursion rule \eqref{E:mrecursion}
and the base cases of the problem. We begin by verifying the base cases:
\begin{align*}
    \det A_{0,1} & = \frac{(2(N+1)p+1)^2}{(2(N+1)p)^2}\left( 1 +\frac{2(N+1)p-1}{2(N+1)p+1}\right) \frac{1}{2(N+1)p+1}= \frac{1}{(N+1)p} ,
\end{align*}
\begin{align*}
    \det A_{0,2} & = \frac{(2(N+1)p+1)^2}{(2(N+1)p-1)^2} \left( 1 +\frac{2(N+1)p-3}{2(N+1)p+1}+ \frac{2(N+1)p-3}{6(N+1)p -1}\right) \frac{1}{2(N+1)p+1} \frac{6(N+1)p-1}{(2(N+1)p)^2}  \\
        &= \frac{28((N+1)p)^2-20(N+1)p-1}{4((N+1)p)^2(2(N+1)p-1)^2} \\
         & = \frac{4}{(N+1)p(2(N+1)p-1)}-\left(\frac{2}{2(N+1)p-1}-\frac{1}{2(N+1)p}\right)^2.
\end{align*}
Now suppose $2 \leq kp+i\leq (N+1)p$. Denote
\begin{align*}
    \alpha_{k,i}
            &=\left(2\lambda_{k,i}-\frac{2 ( \lambda_{k,i}-\lambda_{k,i-1})^2}{ \lambda_{k,i-1}-\lambda_{k,i-2}}+\frac{( \lambda_{k,i}-\lambda_{k,i-1})^2 2\lambda_{k,i-1}}{( \lambda_{k,i-1}-\lambda_{k,i-2})^2} \right) \\
            &=\begin{cases}
                4\frac{(2(N+1)p+1)k-k^2-1}{(2(N+1)p-kp-i)^2},            &kp+i<(N+1)p\\
                3\frac{2((N+1)p)^2+2(N+1)p-1}{(2(N+1)p+1)^2},    &kp+i=(N+1)p
            \end{cases} \\
    \beta_{k,i} &= ( \lambda_{k,i}-\lambda_{k,i-1})^2\left(1-\frac{2\lambda_{k,i-1}}{ \lambda_{k,i-1}-\lambda_{k,i-2}}\right)^2 \\
                &=\begin{cases}
                    \frac{(4 k(N+1)p-2(N+1)p-2k^2+4k-1)^2}{(2(N+1)p-k)^2(2(N+1)p-k+1)^2}, &kp+i<(N+1)p\\
                    \frac{(2(N+1)p^2+2(N+1)p-1)^2}{((N+1)p+1)^2(2(N+1)p+1)^2}, &kp+i=(N+1)p
                \end{cases},
\end{align*}
then the recursion rule \eqref{E:mrecursion} can be written as
\begin{align*}
    & \det A_{k,i}= \alpha_{k,i} \det A_{k,i-1} - \beta_{k,i} \det A_{k,i-2}.
\end{align*}
Further denote
\begin{align*}
    f_{k,i} &= \frac{(2(N+1)p+1)^2}{(2((N+1)p+1)p-kp-i)^2}, \quad kp+i=0,\dots,(N+1)p-1, \\
    g_{k,i} &= 2(N+1)p-2(kp+i)-1, \quad kp+i=0,\dots,(N+1)p-1,\\
    x_{k,i} &= \frac{1}{2(N+1)p+4 (N+1)p (kp+i) - 2 (kp+i)^2+1}, \quad kp+i=0,\dots,(N+1)p-1, \\
    y_{k,i} &= \frac{2(N+1)p+4(N+1)p(kp+i)-2(kp+i)^2+1}{(2(N+1)p+1-kp-i)^2}, \quad kp+i=0,\dots,(N+1)p-1,
\end{align*}
then the solution \eqref{E:mkdet} becomes
\begin{align*}
    \det A_{k,i} =& f_{k,i} \left(1+g_{k,i} \sum_{k'p+i'\leq kp+i} x_{k',i'}\right) \prod_{k'p+i'\leq kp+i} y_{k',i'},
\end{align*}
and \eqref{E:mndet} becomes
\begin{align*}
    \det A_{N,p} =& \frac{(2(N+1)p+1)^2}{((N+1)p+1)^2} \prod_{k'p+i'\leq (N+1)p} y_{k',i'}.
\end{align*}

Substituting \eqref{E:mkdet} in the RHS of \eqref{E:mrecursion} we get that for $kp+i=2,\dots,(N+1)p$
\begin{align*}
     & \alpha_{k,i} \det A_{k,i-1} - \beta_{k,i} \det A_{k,i-2}\\
        &= \alpha_{k,i} f_{k,i-1} \left(1+g_{k,i-1} \sum_{k'p+i'\leq kp+i-1} x_{k',i'}\right) \prod_{k'p+i'\leq kp+i-1} y_{k',i'} \\
        &- \beta_{k,i} f_{k,i-2} \left(1+g_{k,i-2} \sum_{k'p+i'\leq kp+i-2} x_{k',i'}\right) \prod_{k'p+i'\leq kp+i-2} y_{k',i'}\\
        &= \left(\alpha_{k,i} f_{k,i-1}(1+ g_{k,i-1} x_{k,i-1}) -\frac{\beta_{k,i}}{y_{k,i-2}} f_{k,i-2}+\left( \alpha_{k,i} f_{k,i-1} g_{k,i-1}-\frac{\beta_{k,i}}{y_{k,i-1}} f_{k,i-2} g_{k,i-2}\right)  \sum_{k'p+i'\leq kp+i-2} x_{k',i'}\right) \\
        & \times \prod_{k'p+i'\leq kp+i-1} y_{k',i'}.
\end{align*}
It is straightforward (although somewhat involved) to verify that for $kp+i<(N+1)p$
\begin{align*}
    & \alpha_{k,i} f_{k,i-1} (1+ g_{k,i-1} x_{k,i-1}) -\frac{\beta_{k,i-1}}{y_{k,i-1}} f_{k,i-2} = f_{k,i} y_{k,i} (1+g_{k,i} x_{k,i-1}+g_{k,i} x_{k,i}),
\end{align*}
and
\begin{align*}
    & \alpha_{k,i} f_{k,i-1}g_{k,i-1}-\frac{\beta_{k,i}}{y_{k,i-1}} f_{k,i-2} g_{k,i-2} = f_{k,i} g_{k,i} y_{k,i}.
\end{align*}
We therefore get
\begin{align*}
    & \alpha_{k,i} \det A_{k,i-1} -\beta_k \det A_{k,i-2} \\
    & = \left( f_{k,i} y_{k,i}(1+g_{k,i} x_{k,i-1}+g_{k,i} x_{k,i})+ f_{k,i} g_{k,i} y_{k,i} \sum_{k'p+i'\leq kp+i-2} x_{k',i'}\right) \prod_{k'p+i'\leq kp+i-1} y_{k',i'}\\
    %& = f_{k,i} \left( 1+g_{k,i} \sum_{k'p+i'\leq kp+i} x_{k',i'}\right)  \prod_{k'p+i'\leq kp+i} y_{k',i'}\\
    & = \det A_{k,i-1},
\end{align*}
and thus \eqref{E:mkdet} satisfies \eqref{E:mrecursion}. It is also possible to show that
\begin{align*}
    & \alpha_{N,p} f_{N,p-1}(1+ g_{N,p-1} x_{N,p-1}) -\frac{\beta_{N,p}}{y_{N,p-1}} f_{N,p-2} = \frac{(2(N+1)p+1)^2}{((N+1)p+1)^2}, \\
    & \alpha_{N,p} f_{N,p-1} g_{N,p-1}-\frac{\beta_N}{y_{N,p-1}} f_{N,p-2} g_{N,p-2} = 0,
\end{align*}
thus, for $kp+i=(N+1)p$
\begin{align*}
    & \alpha_{N,p}\det A_{N,p-1} -\beta_{N,p} \det A_{{N,p-2}} \\
    & = \frac{(2(N+1)p+1)^2}{((N+1)p+1)^2}  \prod_{k'p+i'\leq (N+1)p+1} y_{k',i'} \\
    & = \det A_{N,p},
\end{align*}
and the expression \eqref{E:mndet} is also verified.

\end{document}

%% file: macros.tex
\newcommand{\BALD}{\begin{aligned}}
\newcommand{\EALD}{\end{aligned}}
\newcommand{\BALDS}{\begin{aligned*}}
\newcommand{\EALDS}{\end{aligned*}}
\newcommand{\BCAS}{\begin{cases}}
\newcommand{\ECAS}{\end{cases}}
\newcommand{\BEAS}{\begin{eqnarray*}}
\newcommand{\EEAS}{\end{eqnarray*}}
\newcommand{\BEQ}{\begin{equation}}
\newcommand{\EEQ}{\end{equation}}
\newcommand{\BIT}{\begin{itemize}}
\newcommand{\EIT}{\end{itemize}}
\newcommand{\BMAT}{\begin{bmatrix}}
\newcommand{\EMAT}{\end{bmatrix}}
\newcommand{\BNUM}{\begin{enumerate}}
\newcommand{\ENUM}{\end{enumerate}}

\newcommand{\BA}{\begin{array}}
\newcommand{\EA}{\end{array}}

\DeclareMathOperator{\trace}{tr}